\documentclass[final,3p]{elsarticle}
 \usepackage{graphics}
 \usepackage{graphicx}
 \usepackage{epsfig}
\usepackage{amssymb}
 \usepackage{amsthm}
 \usepackage{lineno}
 \usepackage{amsmath}
   \numberwithin{equation}{section}
\usepackage{mathrsfs}

\newtheorem{thm}{Theorem}[section]

\newtheorem{lem}[thm]{Lemma}
\newtheorem{prop}[thm]{Proposition}
\newtheorem{defn}[thm]{Definition}

\biboptions{sort&compress,square}
\journal{}
\begin{document}
\begin{frontmatter}
\author[rvt1]{Jian Wang}
\ead{wangj484@nenu.edu.cn}
\author[rvt2]{Yong Wang\corref{cor2}}
\ead{wangy581@nenu.edu.cn}
\cortext[cor2]{Corresponding author.}

\address[rvt1]{School of Science, Tianjin University of Technology and Education, Tianjin, 300222, P.R.China}
\address[rvt2]{School of Mathematics and Statistics, Northeast Normal University,
Changchun, 130024, P.R.China}

\title{  Perturbations of Dirac Operators and A KKW Type Theorem  for \\ Five Dimensional Manifolds with Boundary}

\begin{abstract}
In this paper, we prove a Kastler-Kalau-Walze type theorem associated with  perturbations of Dirac operators
for five dimensional manifolds with boundary.
\end{abstract}
\end{frontmatter}
\section{Introduction}
\label{1}
The noncommutative residue found in \cite{Gu,Wo} plays a prominent role in noncommutative geometry.
 For arbitrary closed compact $n$-dimensional manifolds, the noncommutative reside was introduced by Wodzicki in \cite{Wo} using the theory of zeta
 functions of elliptic pseudodifferential operators.
In \cite{Co1}, Connes used the noncommutative residue to derive a conformal 4-dimensional Polyakov action analogue.
Furthermore, Connes made a challenging observation that the noncommutative residue of the square of the inverse of the
Dirac operator was proportional to the Einstein-Hilbert action in \cite{Co2}.
In \cite{Ka}, Kastler gave a brute-force proof of this theorem. In \cite{KW}, Kalau and Walze proved this theorem in the
normal coordinates system simultaneously, which is called the Kastler-Kalau-Walze theorem now.

An important application of Riemannian geometry is to allow us to define the volume element of a Riemannian manifold $(M_{n} ,g)$.
 The noncommutative residue of Wodzicki \cite{Wo} and Guillemin \cite{Gu} is a trace
on the algebra of (integer order) $\Psi DOs$ on $M$. An important feature is that it
allows us to extend to all $\Psi DOs$ the Dixmier trace, which plays the role of the
integral in the framework of noncommutative geometry. Fedosov etc. defined a noncommutative residue on Boutet de Monvel's algebra
and proved that it was a unique continuous trace in \cite{FGLS}. In \cite{S}, Schrohe gave the relation between the Dixmier trace
and the noncommutative residue for
manifolds with boundary.
For an oriented spin manifold  $M$ with boundary $\partial M$,  by the composition formula in Boutet de Monvel's algebra and
the definition of $ \widetilde{{\rm Wres}}$ \cite{Wa1},  $\widetilde{{\rm Wres}}[(\pi^+D^{-1})^2]$ should be the sum of two terms from
interior and boundary of $M$, where $\pi^+D^{-1}$ is an element in Boutet de Monvel's algebra  \cite{Wa1}.

 Recently, Sitarz and Zajac investigated the spectral action for scalar perturbations of Dirac operators in \cite{SZ}.
Iochum and Levy computed the heat kernel coefficients for Dirac operators with one-form perturbations in \cite{TL}.
In \cite{Wa5}, Wang proved a Kastler-Kalau-Walze type theorem for perturbations of Dirac operators on compact manifolds with or without boundary.
Furthermore, using Dirac operators with perturbations,Wang defined a spectral triple and 
  established an infinitesimal equivariant  index formula in \cite{Wa6,Wa7}.
In \cite{WW1}, we proved a Kastler-Kalau-Walze type theorem for 5-dimensional manifolds with boundary.
Motivated by \cite{SZ,TL,Wa5,Wa6,Wa7,WW1}, we study Dirac operators with one-form perturbations.
In the present paper, we shall restrict our attention to the case of $ \widetilde{{\rm Wres}}[\pi^+(D+c(X))^{-1}\circ \pi^+(D+c(X))^{-1}]$
for five dimensional manifolds with boundary.
Our main result is as follows.

{\bf Main Theorem:} Let $X=X'+a_{n}dx_{n}$ near the boundary and  $X'|_{\partial_{M}}$ is a one form on $\partial_{M}$,
 the following identity holds
\begin{eqnarray}
 \widetilde{{\rm Wres}}[\pi^+(D+c(X))^{-1}\circ \pi^+(D+c(X))^{-1}]
&=&\int_{\partial_{M}}\Big[
\frac{1}{16} \Big( \frac{225}{32}K^{2}+\frac{29}{4}s_{M}\big|_{\partial_{M}}-\big(\frac{155}{12}+5i \big)s_{\partial_{M}}\Big)\nonumber\\
&&-\frac{25}{32} a_{n}K
- |X'|^{2}_{g^{\partial M}}-\frac{35}{64}|X'|^{2}_{g^{\partial M}}K
- 3a_{n}^{2}|_{\partial M}\nonumber\\
&&-\frac{3 }{2} \partial_{ x_{n}}(a_{n})|_{\partial M}
+\frac{15}{8} C_{1}^{1}(\nabla^{\partial M} (X'|_{\partial M})^{*})\Big]\pi^{3}{\rm dvol}_{\partial_{M}},
\end{eqnarray}
where $s_{M}$, $s_{\partial_{M}}$ are respectively scalar curvatures on $M$ and $\partial_{M}$,
 and the vector field $(X'|_{\partial M})^{*}$ is the metric dual of $X'|_{\partial M}$,
  $\nabla^{\partial M}$ is the Levi-civita connection on $\partial M$, $C_{1}^{1}$ is
  the contraction of $(1,1)$ tensors.

 This paper is organized as follows: In Section 2, we define lower dimensional volumes of compact Riemannian manifolds
 with  boundary. In Section 3, for five dimensional spin manifolds with boundary and  the associated Dirac operators   with one-form perturbations,
 we compute $ \widetilde{{\rm Wres}}[\pi^+(D+c(X))^{-1}\circ \pi^+(D+c(X))^{-1}]$ and get a Kastler-Kalau-Walze type theorem in this case.

\section{Lower-Dimensional Volumes of Spin Manifolds  with  boundary}

 In this section we consider an $n$-dimensional oriented Riemannian manifold $(M, g^{M})$ with boundary $\partial_{ M}$ equipped
with a fixed spin structure. We assume that the metric $g^{M}$ on $M$ has
the following form near the boundary
 \begin{equation}
 g^{M}=\frac{1}{h(x_{n})}g^{\partial M}+\texttt{d}x _{n}^{2} ,
\end{equation}
where $g^{\partial M}$ is the metric on $\partial M$. Let $U\subset
M$ be a collar neighborhood of $\partial M$ which is diffeomorphic $\partial M\times [0,1)$. By the definition of $h(x_n)\in C^{\infty}([0,1))$
and $h(x_n)>0$, there exists $\tilde{h}\in C^{\infty}((-\varepsilon,1))$ such that $\tilde{h}|_{[0,1)}=h$ and $\tilde{h}>0$ for some
sufficiently small $\varepsilon>0$. Then there exists a metric $\hat{g}$ on $\hat{M}=M\bigcup_{\partial M}\partial M\times
(-\varepsilon,0]$ which has the form on $U\bigcup_{\partial M}\partial M\times (-\varepsilon,0 ]$
 \begin{equation}
\hat{g}=\frac{1}{\tilde{h}(x_{n})}g^{\partial M}+\texttt{d}x _{n}^{2} ,
\end{equation}
such that $\hat{g}|_{M}=g$.
We fix a metric $\hat{g}$ on the $\hat{M}$ such that $\hat{g}|_{M}=g$.

Let us give the expression of Dirac operators near the boundary. Set  $\widetilde{E}_{n}=\frac{\partial}{\partial x_{n}}$,
 $\widetilde{E}_{j}=\sqrt{h(x_{n})}E_{j}~~(1\leq j \leq n-1)$, where  $\{E_{1},\cdots,E_{n-1}\}$ are orthonormal basis of $T\partial_{M}$.
  Let $\nabla^L$ denote the Levi-civita connection
about $g^M$.
 In the local coordinates $\{x_i; 1\leq i\leq n\}$ and the fixed orthonormal frame $\{\widetilde{E}_{1},\cdots,\widetilde{E}_{n}\}$,
 the connection matrix
 $(\omega_{s,t})$
is defined by
 \begin{equation}
 \nabla^L(\widetilde{E}_{1},\cdots,\widetilde{E}_{n})^{t}= (\omega_{s,t})(\widetilde{E}_{1},\cdots,\widetilde{E}_{n})^{t}.
 \end{equation}

 Usually, the Dirac operator on spin manifold is taken to be the operator, which
comes from the Levi-Civita connection on the tangent bundle to $M$.
So the classical Dirac operator is simply defined on $S$  by
 \begin{equation}
D=\sum_{i}c(e_{i})\nabla _{e_{i}}^{S}
=\sum^n_{j=1}c(\widetilde{E_{j}})\Big[\widetilde{E_j}+\frac{1}{4}\sum_{s,t}\omega_{s,t}(\widetilde{E_j})c(\widetilde{E_s})c(\widetilde{E_t})\Big].
\end{equation}

However, that a wide class of generalized operators which come from connections with  perturbations.
Set $c(X)$ a Clifford action on $M$ and $X$ is a one form,
the Dirac operators with one-form perturbations denote by
  \begin{equation}
\widetilde{D}=\sum_{i}c(e_{i})\nabla _{e_{i}}^{S}+c(X),
\end{equation}
 where $c(e_{i})$ denotes the Clifford action.

The next step is to express the lower dimensional volumes
 $Vol ^{(p_{1},p_{2})}_{n} M$ in a purely differential geometric way.
The local Riemannian invariants makes sense independent of the existence
of a spin structure, we can use the geometric expression for   $Vol ^{(p_{1},p_{2})}_{n} M$
 to extend its definition to general.
Denote by $\mathcal{B}$ Boutet de Monvel's algebra, we recall the main theorem in \cite{FGLS}.
\begin{thm}\label{th:32}{\bf(Fedosov-Golse-Leichtnam-Schrohe)}
 Let $X$ and $\partial X$ be connected, ${\rm dim}X=n\geq3$,
 $A=\left(\begin{array}{lcr}\pi^+P+G &   K \\
T &  S    \end{array}\right)$ $\in \mathcal{B}$ , and denote by $p$, $b$ and $s$ the local symbols of $P,G$ and $S$ respectively.
 Define:
 \begin{eqnarray}
{\rm{\widetilde{Wres}}}(A)&=&\int_X\int_{\bf S}{\rm{tr}}_E\left[p_{-n}(x,\xi)\right]\sigma(\xi)dx \nonumber\\
&&+2\pi\int_ {\partial X}\int_{\bf S'}\left\{{\rm tr}_E\left[({\rm{tr}}b_{-n})(x',\xi')\right]+{\rm{tr}}
_F\left[s_{1-n}(x',\xi')\right]\right\}\sigma(\xi')dx',
\end{eqnarray}
Then~~ a) ${\rm \widetilde{Wres}}([A,B])=0 $, for any
$A,B\in\mathcal{B}$;~~ b) It is a unique continuous trace on
$\mathcal{B}/\mathcal{B}^{-\infty}$.
\end{thm}
 Let $p_{1},p_{2}$ be nonnegative integers and $p_{1}+p_{2}\leq n$. Then by Sec 2.1 of \cite{Wa4},  we have
\begin{defn} Lower-dimensional volumes of spin manifolds with boundary  are defined by
   \begin{equation}\label{}
  {\rm Vol}^{(p_1,p_2)}_n(M):=\widetilde{{\rm Wres}}[\pi^+\widetilde{D}^{-p_1}\circ\pi^+\widetilde{D}^{-p_2}].
\end{equation}
\end{defn}

 Denote by $\sigma_{l}(A)$ the $l$-order symbol of an operator A. An application of (2.1.4) in \cite{Wa1} shows that
\begin{equation}
\widetilde{{\rm Wres}}[\pi^+\widetilde{D}^{-p_1}\circ\pi^+\widetilde{D}^{-p_2}]=\int_M\int_{|\xi|=1}{\rm
trace}_{S(TM)}[\sigma_{-n}(\widetilde{D}^{-p_1-p_2})]\sigma(\xi)\texttt{d}x+\int_{\partial
M}\Phi,
\end{equation}
where
 \begin{eqnarray}
\Phi&=&\int_{|\xi'|=1}\int^{+\infty}_{-\infty}\sum^{\infty}_{j, k=0}
\sum\frac{(-i)^{|\alpha|+j+k+1}}{\alpha!(j+k+1)!}
 {\rm trace}_{S(TM)}
\Big[\partial^j_{x_n}\partial^\alpha_{\xi'}\partial^k_{\xi_n}
\sigma^+_{r}(\widetilde{D}^{-p_1})(x',0,\xi',\xi_n)\nonumber\\
&&\times\partial^\alpha_{x'}\partial^{j+1}_{\xi_n}\partial^k_{x_n}\sigma_{l}
(\widetilde{D}^{-p_2})(x',0,\xi',\xi_n)\Big]d\xi_n\sigma(\xi')\texttt{d}x',
\end{eqnarray}
and the sum is taken over $r-k+|\alpha|+\ell-j-1=-n,r\leq-p_{1},\ell\leq-p_{2}$.

 \section{ A KKW type theorem for five dimensional spin manifolds with boundary }
In this section, we compute the lower dimensional volume for five dimensional compact manifolds with boundary and get a
Kastler-Kalau-Walze type formula in this case.
From now on we always assume that $M$ carries a spin structure so that the spinor bundle and the
Dirac operators with one-form perturbations are defined on $M$.

The following proposition is the key of the computation of lower-dimensional volumes of spin
manifolds with boundary.
\begin{prop}\cite{Wa4}  The following identity holds:
 \begin{eqnarray}
&& 1)~  When~ p_1+p_2=n,~ then, ~ {\rm Vol}^{(p_1,p_2)}_nM=c_0{\rm Vol}_M;\\
&& 2)~  when~  p_1+p_2\equiv n ~ {\rm mod}~  1, ~{\rm Vol}^{(p_1,p_2)}_nM=\int_{\partial M}\Phi.
\end{eqnarray}
\end{prop}
Nextly, for $5$-dimensional spin manifolds with boundary, we compute
${\rm Vol}^{(1,1)}_5$. By Proposition 3.1, we have
 \begin{equation}
\widetilde{{\rm Wres}}[\pi^+(D+c(X))^{-1}\circ \pi^+(D+c(X))^{-1}]=\int_{\partial M}\Phi.
\end{equation}

 Recall  the Dirac operators with one-form perturbations $\widetilde{D}$ of the definitions (2.4) and (2.5).
 Write
  \begin{equation}
D_x^{\alpha}=(-\sqrt{-1})^{|\alpha|}\partial_x^{\alpha};~\sigma(\widetilde{D})=p_1+p_0;
~\sigma(\widetilde{D}^{-1})=\sum^{\infty}_{j=1}q_{-j}.
\end{equation}
By the composition formula of psudodifferential operators, then we have
\begin{eqnarray}
1=\sigma(\widetilde{D}\circ \widetilde{D}^{-1})
&=&\sum_{\alpha}\frac{1}{\alpha!}\partial^{\alpha}_{\xi}[\sigma(\widetilde{D})]D^{\alpha}_{x}[\sigma(\widetilde{D}^{-1})]\nonumber\\
&=&(p_1+p_0)(q_{-1}+q_{-2}+q_{-3}+\cdots)\nonumber\\
& &~~~+\sum_j(\partial_{\xi_j}p_1+\partial_{\xi_j}p_0)(
D_{x_j}q_{-1}+D_{x_j}q_{-2}+D_{x_j}q_{-3}+\cdots)\nonumber\\
&=&p_1q_{-1}+(p_1q_{-2}+p_0q_{-1}+\sum_j\partial_{\xi_j}p_1D_{x_j}q_{-1})\nonumber\\
   &&~~~+(p_1q_{-3}+p_0q_{-2}+\sum_j\partial_{\xi_j}p_1D_{x_j}q_{-2}) +\cdots .
\end{eqnarray}
Thus, we get
\begin{eqnarray}
q_{-1}&=&p_1^{-1};  \\
q_{-2}&=&-p_1^{-1}\Big[p_0p_1^{-1}+\sum_j\partial_{\xi_j}p_1D_{x_j}(p_1^{-1})\Big]; \\
q_{-3}&=&-p_1^{-1}\Big[p_0q_{-2}+\sum_{j=1}^{n-1}c({\rm d}x_{j})\partial_{x_j}q_{-2} +c({\rm d}x_{n})\partial_{x_n}q_{-2}\Big].
\end{eqnarray}

Since $\Phi$ is a global form on $\partial M$, so for any fixed point $x_{0}\in\partial M$, we can choose the normal coordinates
$U$ of $x_{0}$ in $\partial M$(not in $M$) and compute $\Phi(x_{0})$ in the coordinates $\widetilde{U}=U\times [0,1)$ and the metric
$\frac{1}{h(x_{n})}g^{\partial M}+\texttt{d}x _{n}^{2}$.
The dual metric of $g^{M}$ on $\widetilde{U}$ is
$h(x_{n})g^{\partial M}+\frac{\partial}{\partial x_{n}} \otimes\frac{\partial}{\partial x_{n}}.$ Write
$g_{ij}^{M}=g^{M}(\frac{\partial}{\partial x_{i}},\frac{\partial}{\partial x_{j}})$;
$g^{ij}_{M}=g^{M}(d x_{i},dx_{j})$, then

\begin{equation}
[g_{i,j}^{M}]=
\begin{bmatrix}\frac{1}{h( x_{n})}[g_{i,j}^{\partial M}]&0\\0&1\end{bmatrix};\quad
[g^{i,j}_{M}]=\begin{bmatrix} h( x_{n})[g^{i,j}_{\partial M}]&0\\0&1\end{bmatrix},
\end{equation}
and
\begin{equation}
\partial_{x_{s}} g_{ij}^{\partial M}(x_{0})=0,\quad 1\leq i,j\leq n-1;\quad g_{i,j}^{M}(x_{0})=\delta_{ij}.
\end{equation}

Let $\{E_{1},\cdots, E_{n-1}\}$ be an orthonormal frame field in $U$ about $g^{\partial M}$ which is parallel along geodesics and
$E_{i}=\frac{\partial}{\partial x_{i}}(x_{0})$, then $\{\widetilde{E_{1}}=\sqrt{h(x_{n})}E_{1}, \cdots,
\widetilde{E_{n-1}}=\sqrt{h(x_{n})}E_{n-1},\widetilde{E_{n}}=dx_{n}\}$ is the orthonormal frame field in $\widetilde{U}$ about $g^{M}.$
Locally $S(TM)|\widetilde{U}\cong \widetilde{U}\times\wedge^{*}_{C}(\frac{n}{2}).$ Let $\{f_{1},\cdots,f_{n}\}$ be the orthonormal basis of
$\wedge^{*}_{C}(\frac{n}{2})$. Take a spin frame field $\sigma: \widetilde{U}\rightarrow Spin(M)$ such that
$\pi\sigma=\{\widetilde{E_{1}},\cdots, \widetilde{E_{n}}\}$ where $\pi: Spin(M)\rightarrow O(M)$ is a double covering, then
$\{[\sigma, f_{i}], 1\leq i\leq 4\}$ is an orthonormal frame of $S(TM)|_{\widetilde{U}}.$ In the following, since the global form $\Phi$
is independent of the choice of the local frame, so we can compute $\texttt{tr}_{S(TM)}$ in the frame $\{[\sigma, f_{i}], 1\leq i\leq 4\}$.
Let $\{\hat{E}_{1},\cdots,\hat{E}_{n}\}$ be the canonical basis of $R^{n}$ and
$c(\hat{E}_{i})\in Hom(\wedge^{*}_{C}(\frac{n-1}{2}),\wedge^{*}_{C}(\frac{n-1}{2}))$ be the Clifford action. By \cite{Y}, then
\begin{equation}
c(\widetilde{E_{i}})=[(\sigma,c(\hat{E}_{i}))]; \quad c(\widetilde{E_{i}})[(\sigma, f_{i})]=[\sigma,(c(\hat{E}_{i}))f_{i}]; \quad
\frac{\partial}{\partial x_{i}}=[(\sigma,\frac{\partial}{\partial x_{i}})],
\end{equation}
then we have $\frac{\partial}{\partial x_{i}}c(\widetilde{E_{i}})=0$ in the above frame. By Lemma 2.2 in \cite{Wa3}, we have

\begin{lem}\label{le:32}
With the metric $g^{M}$ on $M$ near the boundary
\begin{eqnarray}
\partial_{x_j}(|\xi|_{g^M}^2)(x_0)&=&\left\{
       \begin{array}{c}
        0,  ~~~~~~~~~~ ~~~~~~~~~~ ~~~~~~~~~~~~~{\rm if }~j<n; \\[2pt]
       h'(0)|\xi'|^{2}_{g^{\partial M}},~~~~~~~~~~~~~~~~~~~~~{\rm if }~j=n.
       \end{array}
    \right. \\
\partial_{x_j}[c(\xi)](x_0)&=&\left\{
       \begin{array}{c}
      0,  ~~~~~~~~~~ ~~~~~~~~~~ ~~~~~~~~~~~~~{\rm if }~j<n;\\[2pt]
\partial x_{n}(c(\xi'))(x_{0}), ~~~~~~~~~~~~~~~~~{\rm if }~j=n,
       \end{array}
    \right.
\end{eqnarray}
where $\xi=\xi'+\xi_{n}\texttt{d}x_{n}$.
\end{lem}

 By Lemma 2.1 in \cite{Wa3}, we have
 \begin{lem}\label{le:31}
The symbol of the Dirac operators with one-form perturbations
\begin{eqnarray}
\sigma_{-1}(\widetilde{D}^{-1})&=&q_{-1}=\frac{\sqrt{-1}c(\xi)}{|\xi|^{2}}; \\
\sigma_{-2}(\widetilde{D}^{-1})&=&q_{-2}=\frac{c(\xi)p_{0}c(\xi)}{|\xi|^{4}}+\frac{c(\xi)}{|\xi|^{6}}\sum_{j}c({\rm d} x_{j})
\Big[\partial_{x_{j}}(c(\xi))|\xi|^{2}-c(\xi)\partial_{x_{j}}(|\xi|^{2})\Big]\nonumber\\
&&=\sigma_{-2}(D^{-1})+\frac{c(X)}{|\xi|^{2}}-\frac{2g(X,\xi)c(\xi)}{|\xi|^{4}};\\
\sigma_{-3}(\widetilde{D}^{-1})&=&q_{-3}=-\frac{1}{p_{1}}\Big[p_{0}q_{-2}
+\sum_{j=1}^{n-1}c(dx_{j})\partial_{x_{j}}q_{-2}+c(dx_{n})\partial_{x_{n}}q_{-2}\Big];
\end{eqnarray}
where
 \begin{equation}
p_{0}(x_{0})=-h'(0)c(\texttt{d}x_{n})+c(X).
\end{equation}
\end{lem}

Let us now consider the $q_{-3}$  of the Dirac operators with one-form perturbations. From  Lemma 3.7 in \cite{WW1}, we have
\begin{lem} \cite{WW1} For Dirac operators, the following identity holds:
\begin{eqnarray}
\sigma_{-3}(D^{-1})(x_{0})\Big|_{|\xi'|=1}&=& \frac{-i(h')^{2} }{(1+\xi_{n}^{2})^{3}}c(\xi)c(\texttt{d}x_{n})c(\xi)c(\texttt{d}x_{n})c(\xi) \nonumber\\
&&+\frac{ih' }{(1+\xi_{n}^{2})^{3}}c(\xi)c(\texttt{d}x_{n})c(\xi)c(\texttt{d}x_{n}) \partial_{x_n}[c(\xi')](x_0) \nonumber\\
&&+\frac{-i(h')^{2} }{(1+\xi_{n}^{2})^{4}}c(\xi)c(\texttt{d}x_{n})c(\xi)c(\texttt{d}x_{n})c(\xi) \nonumber\\
&&+\frac{1}{8}\sum_{ \beta i s \alpha}R^{\partial_{M}}_{ \beta i s \alpha}(x_0)
 \frac{ -i}{(1+\xi_{n}^{2})^{3}}c(\xi)c(\widetilde{E}_{i})c(\xi)c(\widetilde{E}_{\beta})c(\widetilde{E}_{s})c(\widetilde{E}_{\alpha})c(\xi)\nonumber\\
&&+\frac{1}{6} \sum_{l,t<n}\xi_{l}\Big(R^{\partial_{M}}_{tilj}(x_0)+R^{\partial_{M}}_{tjli}(x_0)\Big)
 \frac{-i }{(1+\xi_{n}^{2})^{3}}c(\xi)c(\widetilde{E}_{i})c(\xi)c(\texttt{d}x_{j})c(\widetilde{E}_{t}) \nonumber\\
 &&+\frac{1}{3} \sum_{\alpha ,\beta <n}\Big(R^{\partial_{M}}_{i\alpha j\beta}(x_0)+R^{\partial_{M}}_{i\beta j\alpha}(x_0)\Big)
 \xi_{\alpha}\xi_{\beta}\frac{i }{(1+\xi_{n}^{2})^{4}}c(\xi)c(\widetilde{E}_{i}) c(\xi)c(\texttt{d}x_{j})c(\xi)\nonumber\\
 &&+i\Big( \frac{h' }{(1+\xi_{n}^{2})^{3}}+ \frac{h' }{(1+\xi_{n}^{2})^{4}}\Big)
    c(\xi)c(\texttt{d}x_{n})\partial_{x_n}[c(\xi')](x_0)c(\texttt{d}x_{n})c(\xi)\nonumber\\
   &&-i \Big( \frac{(h')^{2}-h'' }{(1+\xi_{n}^{2})^{3}}+\frac{2(h')^{2}-h'' }{(1+\xi_{n}^{2})^{4}}+\frac{3(h')^{2} }{(1+\xi_{n}^{2})^{5}}\Big)
    c(\xi)c(\texttt{d}x_{n})c(\xi)c(\texttt{d}x_{n})c(\xi)\nonumber\\
   &&+ i\Big( \frac{h' }{(1+\xi_{n}^{2})^{3}}+ \frac{3h' }{(1+\xi_{n}^{2})^{4}}\Big)
    c(\xi)c(\texttt{d}x_{n})c(\xi)c(\texttt{d}x_{n})\partial_{x_n}[c(\xi')](x_0)\nonumber\\
   &&+\frac{-i}{(1+\xi_{n}^{2})^{3}}c(\xi)c(\texttt{d}x_{n})\partial_{x_n}[c(\xi')](x_0)  c(\texttt{d}x_{n}) \partial_{x_n}[c(\xi')](x_0)\nonumber\\
   &&+\Big(\frac{3}{4}(h'(0))^{2}-\frac{1}{2}h''(0)\Big)\frac{-i }{(1+\xi_{n}^{2})^{4}} c(\xi)c(\texttt{d}x_{n})c(\xi)c(\texttt{d}x_{n})c(\xi').
\end{eqnarray}
\end{lem}

 From (3.16) and  Lemma 3.2-Lemma 3.4, we obtain
 \begin{lem} For Dirac operators with one-form perturbations, the following identity holds:
\begin{eqnarray}
\sigma_{-3}(\widetilde{D}^{-1})(x_{0})\Big|_{|\xi'|=1}&=& \sigma_{-3}({D}^{-1})(x_{0})\Big|_{|\xi'|=1}
+\Big[-q_{1}c(X)\sigma_{-2}({D}^{-1})-q_{1}\frac{c(X)}{|\xi|^{2}}+q_{1}\frac{2g(X,\xi)c(\xi)}{|\xi|^{4}}\nonumber\\
&&-q_{1}\sum_{j=1}^{n-1}c(dx_{j})\partial_{x_{j}}(\frac{c(X)}{|\xi|^{2}})
+q_{1}\sum_{j=1}^{n-1}c(dx_{j})\partial_{x_{j}}(\frac{2g(X,\xi)c(\xi)}{|\xi|^{4}})\nonumber\\
&&-q_{1}c(dx_{n})\partial_{x_{n}}(\frac{c(X)}{|\xi|^{2}})+q_{1}\partial_{x_{n}}(\frac{2g(X,\xi)c(\xi)}{|\xi|^{4}})\Big](x_{0})\Big|_{|\xi'|=1}\nonumber\\
&:=&\sigma_{-3}({D}^{-1})(x_{0})\Big|_{|\xi'|=1}+R_{-3}(x_{0})\Big|_{|\xi'|=1}.
\end{eqnarray}
\end{lem}

From the remark above, now we can compute $\Phi$ (see the formula (2.9) for the definition of $\Phi$).
Since the sum is taken over $-r-\ell+k+j+|\alpha|-1=5, \ r, \ell\leq-1$, then we have the $\int_{\partial_{M}}\Phi$ is the sum of the following
 fifteen cases. Such cases (1)-(6) have been studied,
similar to \textbf{Case (1)}-\textbf{Case (6)} in \cite{WW1}, we obtain

\textbf{Case (1)}: \ $r=-1, \ \ell=-1, \ k=0, \ j=1, \ |\alpha|=1$

From (2.9), we have
\begin{equation}
\text{ Case \ (1)}=\frac{i}{2}\int_{|\xi'|=1}\int_{-\infty}^{+\infty}\sum_{|\alpha|=1}\text{trace}
\Big[\partial_{x_{n}}\partial_{\xi'}^{\alpha}\pi_{\xi_{n}}^{+}q_{-1}\times
\partial_{x'}^{\alpha}\partial_{\xi_{n}}^{2}q_{-1}\Big](x_{0})\texttt{d}\xi_{n}\sigma(\xi')\texttt{d}x' .
\end{equation}

By Lemma 3.3, for $i<n$, we have
\begin{equation}
\partial_{x_i}q_{-1}(x_0)=\partial_{x_i}\left(\frac{\sqrt{-1}c(\xi)}{|\xi|^2}\right)(x_0)=
\frac{\sqrt{-1}\partial_{x_i}[c(\xi)](x_0)}{|\xi|^2}
-\frac{\sqrt{-1}c(\xi)\partial_{x_i}(|\xi|^2)(x_0)}{|\xi|^4}=0.
\end{equation}
So Case (1) vanishes.

\textbf{Case (2)}: \  $r=-1, \ \ell=-1, \ k=0, \ j=2, \ |\alpha|=0$

From (2.9), we have
\begin{equation}
\text{ Case \ (2)}=\frac{i}{6}\int_{|\xi'|=1}\int_{-\infty}^{+\infty}\sum_{j=2}\text{trace}\Big[\partial_{x_{n}}^{2}\pi_{\xi_{n}}^{+}
q_{-1}\times\partial_{\xi_{n}}^{3}q_{-1}\Big](x_{0})\texttt{d}\xi_{n}\sigma(\xi')\texttt{d}x'.
\end{equation}
By Lemma 3.2, a simple computation shows
\begin{equation}
\partial_{\xi_{n}}^{3}q_{-1}(x_{0})\Big|_{|\xi'|=1}
=\frac{24\xi_{n}-24\xi_{n}^{3}}{(1+\xi_{n}^{2})^{4}}\sqrt{-1}c(\xi')+\frac{-6\xi_{n}^{4}
+36\xi_{n}^{2}-6}{(1+\xi_{n}^{2})^{4}}\sqrt{-1}c(\texttt{d}x_{n}).
\end{equation}
From Lemma 3.2, Lemma 3.3 and Lemma 3.4, we obtain
\begin{eqnarray}
\partial_{x_{n}}^{2}\pi_{\xi_{n}}^{+}q_{-1}(x_{0})\Big|_{|\xi'|=1}
&=&\Big(\frac{3}{4}(h'(0))^{2}-\frac{1}{2}h''(0)\Big)\frac{c(\xi')}{2(\xi_n-i)}
-h'(0) \frac{\xi_n -2i }{2(\xi_n-i)^{2}}\partial_{x_{n}}\big(c(\xi')\big) \nonumber\\
&&-h''(0)\Big[\frac{\xi_n -2i }{4(\xi_n-i)^{2}}c(\xi')+\frac{1 }{4(\xi_n-i)^2}c(dx_n)  \Big] \nonumber\\
&&+2i(h'(0))^{2}\Big[\frac{-3i\xi_n^{2} -9\xi_n+8i}{16(\xi_n-i)^{3}}c(\xi')
+\frac{-i\xi_n-3}{16(\xi_n-i)^{3}}c(dx_n)  \Big].
\end{eqnarray}
Using the Clifford relations combined with the cyclicity of the trace
and $\texttt{tr}AB=\texttt{tr}BA$, then
\begin{eqnarray}
&&\texttt{tr}[c(\xi')c(\texttt{d}x_{n})]=0; \ \texttt{tr}[c(\texttt{d}x_{n})^{2}]=-4;\ \texttt{tr}[c(\xi')^{2}](x_{0})|_{|\xi'|=1}=-4;\nonumber\\
&&\texttt{tr}[\partial_{x_{n}}[c(\xi')]c(\texttt{d}x_{n})]=0; \ \texttt{tr}[\partial_{x_{n}}c(\xi')\times c(\xi')](x_{0})|_{|\xi'|=1}=-2h'(0).
\end{eqnarray}
From (3.22)-(3.25) and direct computations, we obtain
\begin{eqnarray}
&&\text{trace}\Big[\partial_{x_{n}}^{2}\pi_{\xi_{n}}^{+}
q_{-1}\times\partial_{\xi_{n}}^{3}q_{-1}\Big](x_{0})\Big|_{|\xi'|=1}\nonumber\\
&=&i \big(h'(0)\big)^{2}\frac{3(33\xi_{n}^{5}-75i\xi_{n}^{4}-94\xi_{n}^{3}+90i\xi_{n}^{2}+57\xi_{n}-3i)}
 {2(\xi_{n}-i)^{3}(1+\xi_{n}^{2})^{4}}\nonumber\\
 &&+i h''(0)\frac{6(-9\xi_{n}^{4}+12i\xi_{n}^{3}+14\xi_{n}^{2}-12i\xi_{n}-1)}
 {2(\xi_{n}-i)^{2}(1+\xi_{n}^{2})^{4}}.
\end{eqnarray}
Therefore
\begin{eqnarray}
\text{ Case \ (2) }
&=&-\frac{1}{6}\big(h'(0)\big)^{2}\int_{|\xi'|=1}\int_{-\infty}^{+\infty}\frac{3(33\xi_{n}^{5}-75i\xi_{n}^{4}-94\xi_{n}^{3}
+90i\xi_{n}^{2}+57\xi_{n}-3i)} {2(\xi_{n}-i)^{3}(1+\xi_{n}^{2})^{4}}\texttt{d}\xi_{n}\sigma(\xi')\texttt{d}x'\nonumber\\
&&-\frac{1}{6}h''(0)\int_{|\xi'|=1}\int_{-\infty}^{+\infty}\frac{6(-9\xi_{n}^{4}+12i\xi_{n}^{3}+14\xi_{n}^{2}-12i\xi_{n}-1)}
 {2(\xi_{n}-i)^{2}(1+\xi_{n}^{2})^{4}}\texttt{d}\xi_{n}\sigma(\xi')\texttt{d}x'\nonumber\\
 &=&-\frac{1}{6}\big(h'(0)\big)^{2} \Omega_{3}\int_{\Gamma^{+}}\frac{3(33\xi_{n}^{5}-75i\xi_{n}^{4}-94\xi_{n}^{3}
+90i\xi_{n}^{2}+57\xi_{n}-3i)} {2(\xi_{n}-i)^{3}(1+\xi_{n}^{2})^{4}}\texttt{d}\xi_{n}\texttt{d}x'\nonumber\\
&&-\frac{1}{6}h''(0)\Omega_{3}\int_{\Gamma^{+}}\frac{6(-9\xi_{n}^{4}+12i\xi_{n}^{3}+14\xi_{n}^{2}-12i\xi_{n}-1)}
 {2(\xi_{n}-i)^{2}(1+\xi_{n}^{2})^{4}}\texttt{d}\xi_{n}\texttt{d}x'\nonumber\\
&=&-\frac{1}{6}\big(h'(0)\big)^{2}  \frac{ 2\pi i }{6!}\bigg[\frac{3(33\xi_{n}^{5}-75i\xi_{n}^{4}-94\xi_{n}^{3}
+90i\xi_{n}^{2}+57\xi_{n}-3i)}{2(\xi_{n}+i)^{4}}\bigg]^{(6)}\bigg|_{\xi_{n}=i}\Omega_{3}\texttt{d}x'  \nonumber\\
&&-\frac{1}{6}h''(0) \frac{ 2\pi i }{5!}\bigg[\frac{6(-9\xi_{n}^{4}+12i\xi_{n}^{3}+14\xi_{n}^{2}-12i\xi_{n}-1)}
{(\xi_{n}+i)^{4}}\bigg]^{(5)}\bigg|_{\xi_{n}=i}\Omega_{3}\texttt{d}x'  \nonumber\\
&=&\Big(\frac{29}{64}\big(h'(0)\big)^{2}-\frac{3}{8}h''(0)\Big)\pi\Omega_{3}\texttt{d}x',
\end{eqnarray}
where $\Omega_{3}$ is the canonical volume of $S^{3}$.

Similarly, we have

\textbf{Case (3)}: \ $r=-1, \ \ell=-1, \ k=0, \ j=0, \ |\alpha|=2$

\begin{eqnarray}
\text{ Case \ (3)}&=&\frac{i}{2}\int_{|\xi'|=1}\int_{-\infty}^{+\infty}\sum_{|\alpha|=2}
\text{trace}\Big[\partial_{\xi'}^{\alpha}\pi_{\xi_{n}}^{+}q_{-1}\times
\partial_{x'}^{\alpha}\partial_{\xi_{n}}q_{-1}\Big](x_{0})\texttt{d}\xi_{n}\sigma(\xi')\texttt{d}x'\nonumber\\
&=&-\frac{1}{4}s_{\partial_{M}}\pi^{3}\texttt{d}x',
\end{eqnarray}
where $\sum_{t,l <n}R^{\partial_{M}}_{tltl}(x_0) $ is the scalar curvature $s_{\partial_{M}}$.

\textbf{Case (4)}: \ $r=-1, \ \ell=-1, \ k=1, \ j=1, \ |\alpha|=0$
\begin{eqnarray}
\text{ Case \ (4)}&=&\frac{i}{6}\int_{|\xi'|=1}\int_{-\infty}^{+\infty}
\text{trace}\Big[\partial_{x_{n}}\partial_{\xi_{n}}\pi_{\xi_{n}}^{+}q_{-1}\times
\partial_{\xi_{n}}^{2}\partial_{x_{n}}q_{-1}\Big](x_{0})\texttt{d}\xi_{n}\sigma(\xi')\texttt{d}x' \nonumber\\
  &=&-\frac{i}{6}\int_{|\xi'|=1}\int_{-\infty}^{+\infty}
\text{trace}\Big[\partial_{x_{n}}\pi_{\xi_{n}}^{+}q_{-1}\times
\partial_{\xi_{n}}^{3}\partial_{x_{n}}q_{-1}\Big](x_{0})\texttt{d}\xi_{n}\sigma(\xi')\texttt{d}x'\nonumber\\
&=&-\frac{5}{16}\big(h'(0)\big)^{2}\pi\Omega_{3}\texttt{d}x'.
\end{eqnarray}

\textbf{Case (5)}: \ $r=-1, \ \ell=-1, \ k=1, \ j=0, \ |\alpha|=1$

\begin{eqnarray}
\text{ Case \ (5)}&=&\frac{i}{2}\int_{|\xi'|=1}\int_{-\infty}^{+\infty}\sum_{|\alpha|=1}
\text{trace}\Big[\partial_{\xi'}^{\alpha}\partial_{\xi_{n}}\pi_{\xi_{n}}^{+}q_{-1}\times
\partial_{x'}^{\alpha}\partial_{\xi_{n}}\partial_{x_{n}}q_{-1}\Big](x_{0})\texttt{d}\xi_{n}\sigma(\xi')\texttt{d}x' \nonumber\\
&=&0.
\end{eqnarray}

\textbf{Case (6)}: \ $r=-1, \ \ell=-1, \ k=2, \ j=0, \ |\alpha|=0$

\begin{eqnarray}
\text{ Case \ (6)}&=&\frac{i}{6}\int_{|\xi'|=1}\int_{-\infty}^{+\infty}\sum_{k=2}
\text{trace}\Big[\partial_{\xi_{n}}^{2}\pi_{\xi_{n}}^{+}q_{-1}
\partial_{\xi_{n}}\partial_{x_{n}}^{2}q_{-1}\Big](x_{0})\texttt{d}\xi_{n}\sigma(\xi')\texttt{d}x' \nonumber\\
&=&\Big(\frac{29}{64}\big(h'(0)\big)^{2}-\frac{3}{8}h''(0)\Big)\pi\Omega_{3}\texttt{d}x'.
\end{eqnarray}

Now we discuss the cases (7)-(15).

\textbf{Case (7)}: \ $r=-1, \ \ell=-2, \ k=0, \ j=1, \ |\alpha|=0$

From (2.9) and the Leibniz rule, we obtain
\begin{eqnarray}
\text{ Case \ (7)}&=&\frac{1}{2}\int_{|\xi'|=1}\int_{-\infty}^{+\infty}
\text{trace}\Big[\partial_{\xi_{n}}\partial_{x_{n}}\pi_{\xi_{n}}^{+}q_{-1}\times
\partial_{\xi_{n}}q_{-2}\Big](x_{0})\texttt{d}\xi_{n}\sigma(\xi')\texttt{d}x' \nonumber\\
  &=&-\frac{1}{2}\int_{|\xi'|=1}\int_{-\infty}^{+\infty}
\text{trace}\Big[\partial_{\xi_{n}}^{2}\partial_{x_{n}}\pi_{\xi_{n}}^{+}q_{-1}
\times q_{-2}\Big](x_{0})\texttt{d}\xi_{n}\sigma(\xi')\texttt{d}x'.
\end{eqnarray}
By Lemma 3.3,(2.2.22) in \cite{Wa3} and direct computations, we obtain
\begin{equation}
\partial_{\xi_{n}}^{2}\partial_{x_{n}}\pi_{\xi_{n}}^{+}q_{-1}(x_0)|_{|\xi'|=1}
=\frac{1}{(\xi_n-1)^{3}}\partial_{x_n}[c(\xi')](x_0)
+h'(0)\Big[\frac{4i-\xi_n}{2(\xi_n-1)^{4}}c(\xi')
-\frac{3}{2(\xi_n-1)^{4}}c(\texttt{d}x_n)\Big].
\end{equation}
From Lemma 3.2 and Lemma 3.3, we have
\begin{eqnarray}
q_{-2}&=&\frac{c(\xi)p_{0}c(\xi)}{|\xi|^{4}}+\frac{c(\xi)}{|\xi|^{6}}\sum_{j}c({\rm d} x_{j})
\Big[\partial_{x_{j}}(c(\xi))|\xi|^{2}-c(\xi)\partial_{x_{j}}(|\xi|^{2})\Big]\nonumber\\
&&=\sigma_{-2}(D^{-1})+\frac{c(X)}{|\xi|^{2}}-\frac{2g(X,\xi)c(\xi)}{|\xi|^{4}}.
\end{eqnarray}
Then
\begin{eqnarray}
\text{ Case \ (7)}&=&-\frac{1}{2}\int_{|\xi'|=1}\int_{-\infty}^{+\infty}
\text{trace}\Big[\partial_{\xi_{n}}^{2}\partial_{x_{n}}\pi_{\xi_{n}}^{+}q_{-1}
\times \sigma_{-2}(D^{-1})\Big](x_{0})\texttt{d}\xi_{n}\sigma(\xi')\texttt{d}x'\nonumber\\
&&-\frac{1}{2}\int_{|\xi'|=1}\int_{-\infty}^{+\infty}
\text{trace}\Big[\partial_{\xi_{n}}^{2}\partial_{x_{n}}\pi_{\xi_{n}}^{+}q_{-1}
\times \Big(\frac{c(X)}{|\xi|^{2}}-\frac{2g(X,\xi)c(\xi)}{|\xi|^{4}}\Big)\Big](x_{0})\texttt{d}\xi_{n}\sigma(\xi').
\end{eqnarray}
Let $X=X'+a_{n}dx_{n}$, by the relation of the Clifford action and ${\rm tr}{AB}={\rm tr }{BA}$, then we have the equalities:
\begin{eqnarray}
&&{\rm tr}[c(\xi')c(X)]=-4g(X,\xi');~~{\rm tr}[c(dx_n)c(X)]=-4a_{n}. \nonumber
\end{eqnarray}
Considering for $i<n$,  $\int_{|\xi'|=1}\{\xi_{i_1}\xi_{i_2}\cdots\xi_{i_{2d+1}}\}\sigma(\xi')=0$.
 From (3.34), (3.35) and direct computations, we obtain
\begin{equation}
\text{trace}\Big[\partial_{\xi_{n}}^{2}\partial_{x_{n}}\pi_{\xi_{n}}^{+}q_{-1}
\times \Big(\frac{c(X)}{|\xi|^{2}}-\frac{2g(X,\xi)c(\xi)}{|\xi|^{4}}\Big)\Big](x_{0})
=a_{n}h'(0)\frac{-6\xi_{n}^{2}+8i\xi_{n}+6} {(\xi_{n}-i)^{6}(\xi_{n}+i)^{2}}.
\end{equation}
By \text{ Case \ (7) } in \cite{WW1}, then
\begin{eqnarray}
\text{ Case \ (7) }
&=&\frac{39}{32}\big(h'(0)\big)^{2}\pi\Omega_{3}\texttt{d}x'\nonumber\\
&&-\frac{1}{2}\int_{|\xi'|=1}\int_{-\infty}^{+\infty}
\text{trace}\Big[\partial_{\xi_{n}}^{2}\partial_{x_{n}}\pi_{\xi_{n}}^{+}q_{-1}
\times \Big(\frac{c(X)}{|\xi|^{2}}-\frac{2g(X,\xi)c(\xi)}{|\xi|^{4}}\Big)\Big](x_{0})\texttt{d}\xi_{n}\sigma(\xi')\nonumber\\
 &=&\frac{39}{32}\big(h'(0)\big)^{2}\pi\Omega_{3}\texttt{d}x'-\frac{1}{2}a_{n}h'(0)
\frac{ 2\pi i }{5!}\bigg[\frac{-6\xi_{n}^{2}+8i\xi_{n}+6}{(\xi_{n}+i)^{2}}\bigg]^{(5)}\bigg|_{\xi_{n}=i}\Omega_{3}\texttt{d}x' \nonumber\\
&=&\frac{39}{32}\big(h'(0)\big)^{2}\pi\Omega_{3}\texttt{d}x'-\frac{5}{8}a_{n}h'(0)\pi\Omega_{3}\texttt{d}x'.
\end{eqnarray}

\textbf{Case (8)}: \ $r=-1, \ \ell=-2, \ k=0, \ j=0, \ |\alpha|=1$

From (2.9) and the Leibniz rule, we obtain
\begin{eqnarray}
\text{ Case \ (8)}&=&-\int_{|\xi'|=1}\int_{-\infty}^{+\infty}\sum_{|\alpha|=1}
\text{trace}\Big[\partial_{\xi'}^{\alpha}\pi_{\xi_{n}}^{+}q_{-1}\times
\partial_{x'}^{\alpha}\partial_{\xi_{n}}q_{-2}\Big](x_{0})\texttt{d}\xi_{n}\sigma(\xi')\texttt{d}x' \nonumber\\
  &=&\int_{|\xi'|=1}\int_{-\infty}^{+\infty}\sum_{|\alpha|=1}
\text{trace}\Big[\partial_{\xi_{n}}\partial_{\xi'}^{\alpha}\pi_{\xi_{n}}^{+}q_{-1}\times
\partial_{x'}^{\alpha}q_{-2}\Big](x_{0})\texttt{d}\xi_{n}\sigma(\xi')\texttt{d}x'\nonumber\\
  &=&\int_{|\xi'|=1}\int_{-\infty}^{+\infty}\sum_{|\alpha|=1}
\text{trace}\Big[\partial_{\xi_{n}}\partial_{\xi'}^{\alpha}\pi_{\xi_{n}}^{+}q_{-1}\times
\partial_{x'}^{\alpha}\sigma_{-2}(D^{-1})\Big](x_{0})\texttt{d}\xi_{n}\sigma(\xi')\texttt{d}x'\nonumber\\
&&+\int_{|\xi'|=1}\int_{-\infty}^{+\infty}\sum_{|\alpha|=1}
\text{trace}\Big[\partial_{\xi_{n}}\partial_{\xi'}^{\alpha}\pi_{\xi_{n}}^{+}q_{-1}\times
\partial_{x'}^{\alpha}\Big(\frac{c(X)}{|\xi|^{2}}-\frac{2g(X,\xi)c(\xi)}{|\xi|^{4}}\Big)\Big](x_{0})
\texttt{d}\xi_{n}\sigma(\xi')\texttt{d}x'.\nonumber\\
\end{eqnarray}
From Lemma 3.2 and Lemma 3.3, a simple computation shows
\begin{equation}
\partial_{\xi_{n}}\pi_{\xi_{n}}^{+}\partial_{\xi'}^{\alpha}q_{-1}(x_{0})\Big|_{|\xi'|=1}
=\frac{-1}{2(\xi_{n}-i)^{2}}c(\texttt{d}x_{i})-\xi_{i}\frac{3i-\xi_{n}}{2(\xi_{n}-i)^{3}}c(\xi')
+\xi_{i}\frac{1}{(\xi_{n}-i)^{3}}c(\texttt{d}x_{n}).
\end{equation}
From Lemma 3.2, we have
\begin{equation}
\partial_{x'}^{\alpha}\Big(\frac{c(X)}{|\xi|^{2}}-\frac{2g(X,\xi)c(\xi)}{|\xi|^{4}}\Big)\Big|_{|\xi'|=1}
=\frac{1}{|\xi|^{2}}\partial_{x'}[c(X)]-
\frac{2c(\xi)}{|\xi|^{4}}\partial_{x'}[g(X,\xi)].
\end{equation}
By the relation of the Clifford action and ${\rm tr}{AB}={\rm tr }{BA}$, then we have the equalities at a fixed point $x_{0}$:
\begin{eqnarray}
&&{\rm tr}\Big[c(\xi')\partial_{x'}[c(X)]\Big](x_{0})=-4\partial_{x'}[g(X,\xi')](x_{0}). \nonumber
\end{eqnarray}
Combining (3.39), (3.40) and direct computations, we obtain
\begin{eqnarray}
&&\sum_{|\alpha|=1}
\text{trace}\Big[\partial_{\xi_{n}}\partial_{\xi'}^{\alpha}\pi_{\xi_{n}}^{+}q_{-1}\times
\partial_{x'}^{\alpha}\Big(\frac{c(X)}{|\xi|^{2}}-\frac{2g(X,\xi)c(\xi)}{|\xi|^{4}}\Big)\Big](x_{0})\nonumber\\
&=&\frac{2} {(\xi_{n}-i)^{3}(\xi_{n}+i)}\sum_{j=1}^{n-1}\partial_{x_{j}}[g(X,dx_{j})]
+\frac{-2\xi_{n}^{3}+6i\xi_{n}^{2}-2\xi_{n}-2i} {(\xi_{n}-i)^{5}(\xi_{n}+i)^{2}}\sum_{i,j=1}^{n-1}\xi_{i}\xi_{j}\partial_{x_{i}}[g(X,dx_{j})].
\end{eqnarray}
Then an application of (16) in \cite{Ka} shows
\begin{equation}
\int_{S^{3}} \xi_{\mu}\xi_{\nu}=\frac{\pi^{2}}{2}\delta^{\mu\nu}.
\end{equation}
Using the integration over $S^{3}$ and the shorthand $ \int=\frac{1}{2\pi^{2}}\int_{S^{3}}d^{3}\nu$, we obtain  $\Omega_{3}=2\pi^{2} $.
Considering for $i<n$,  $\int_{|\xi'|=1}\{\xi_{i_1}\xi_{i_2}\cdots\xi_{i_{2d+1}}\}\sigma(\xi')=0$.
Let $X=X'+a_{n}dx_{n}$ near the boundary and  $X'$ is a one form on $\partial_{M}$, therefore
\begin{eqnarray}
\text{ Case \ (8) }
&=&\Big(\frac{3 }{16}-\frac{5}{32} i\Big)s_{\partial_{M}}\pi\Omega_{3}\texttt{d}x'\nonumber\\
&&+\int_{|\xi'|=1}\int_{-\infty}^{+\infty}\sum_{|\alpha|=1}
\text{trace}\Big[\partial_{\xi_{n}}\partial_{\xi'}^{\alpha}\pi_{\xi_{n}}^{+}q_{-1}\times
\partial_{x'}^{\alpha}\Big(\frac{c(X)}{|\xi|^{2}}-\frac{2g(X,\xi)c(\xi)}{|\xi|^{4}}\Big)\Big](x_{0})
\texttt{d}\xi_{n}\sigma(\xi')\texttt{d}x'\nonumber\\
&&=\Big(\frac{3 }{16}-\frac{5}{32} i\Big)s_{\partial_{M}}\pi\Omega_{3}\texttt{d}x'
-\frac{9}{16}\pi\Omega_{3} \sum_{j=1}^{n-1}\partial_{x_{j}}[g(X,dx_{j})]\texttt{d}x'\nonumber\\
&&=\Big(\frac{3 }{16}-\frac{5}{32} i\Big)s_{\partial_{M}}\pi\Omega_{3}\texttt{d}x'
-\frac{9}{16}\pi\Omega_{3}C_{1}^{1}(\nabla^{\partial M} (X'|_{\partial M})^{*}))\texttt{d}x',
\end{eqnarray}
where the vector field $(X'|_{\partial M})^{*}=g^{\partial M}(X'|_{\partial M}),\cdot)$ and
  $\nabla^{\partial M}$ is the Levi-civita connection on $\partial M$, $C_{1}^{1}$ is
  the contraction of $(1,1)$ tensors.

\textbf{Case (9)}: \ $r=-1, \ \ell=-2, \ k=1, \ j=0, \ |\alpha|=0$

From (2.9) and the Leibniz rule, we obtain
\begin{eqnarray}
\text{ Case \ (9)}&=&-\frac{1}{2}\int_{|\xi'|=1}\int_{-\infty}^{+\infty}\sum_{|\alpha|=1}
\text{trace}\Big[\partial_{\xi_{n}}\pi_{\xi_{n}}^{+}q_{-1}\times
\partial_{\xi_{n}}\partial_{x_{n}}q_{-2})\Big](x_{0})\texttt{d}\xi_{n}\sigma(\xi')\texttt{d}x' \nonumber\\
  &=&\frac{1}{2}\int_{|\xi'|=1}\int_{-\infty}^{+\infty}\sum_{|\alpha|=1}
\text{trace}\Big[\partial_{\xi_{n}}^{2}\pi_{\xi_{n}}^{+}q_{-1}\times
\partial_{x_{n}}q_{-2}\Big](x_{0})\texttt{d}\xi_{n}\sigma(\xi')\texttt{d}x'\nonumber\\
&=&\frac{1}{2}\int_{|\xi'|=1}\int_{-\infty}^{+\infty}\sum_{|\alpha|=1}
\text{trace}\Big[\partial_{\xi_{n}}^{2}\pi_{\xi_{n}}^{+}q_{-1}\times
\partial_{x_{n}}\sigma_{-2}(D^{-1})\Big](x_{0})\texttt{d}\xi_{n}\sigma(\xi')\texttt{d}x'\nonumber\\
&&+\frac{1}{2}\int_{|\xi'|=1}\int_{-\infty}^{+\infty}\sum_{|\alpha|=1}
\text{trace}\Big[\partial_{\xi_{n}}^{2}\pi_{\xi_{n}}^{+}q_{-1}\times
\partial_{x_{n}}\Big(\frac{c(X)}{|\xi|^{2}}-\frac{2g(X,\xi)c(\xi)}{|\xi|^{4}}\Big)
\Big](x_{0})\texttt{d}\xi_{n}\sigma(\xi')\texttt{d}x'.\nonumber\\
\end{eqnarray}
By (2.2.29) in \cite{Wa3}, we have
\begin{equation}
\partial_{\xi_n}^{2}\pi^+_{\xi_n}q_{-1}(x_0)|_{|\xi'|=1}=\frac{1}{(\xi_n-i)^{3}}c(\xi')+\frac{i}{2(\xi_n-i)^{3}}c(\texttt{d}x_n).
\end{equation}
From Lemma 3.2, a simple computation shows
\begin{eqnarray}
\partial_{x_{n}}\Big(\frac{c(X)}{|\xi|^{2}}-\frac{2g(X,\xi)c(\xi)}{|\xi|^{4}}\Big)
&=&\partial_{x_{n}}\Big(\frac{c(X)}{|\xi|^{2}}\Big)
-\partial_{x_{n}}\Big(\frac{2g(X,\xi)c(\xi)}{|\xi|^{4}}\Big)\nonumber\\
&=&\frac{1}{|\xi|^{2}}\partial_{x_{n}}\Big(c(X)\Big)-\frac{h'(0)|\xi'|^{2}}{|\xi|^{4}}c(X)
-\frac{2c(\xi)}{|\xi|^{4}}\partial_{x_{n}}\Big(g(X,\xi)\Big)\nonumber\\
&&-\frac{2g(X,\xi)}{|\xi|^{4}}\partial_{x_{n}}\Big(c(\xi)\Big)
+\frac{4g(X,\xi)c(\xi)}{|\xi|^{6}}\partial_{x_{n}}\Big(|\xi|^{2})\Big).
\end{eqnarray}
By the relation of the Clifford action and ${\rm tr}{AB}={\rm tr }{BA}$, then we have the equalities:
\begin{eqnarray}
&&{\rm tr}\Big[c(dx_{n})\partial_{x_{n}}[c(X)]\Big]=-4\partial_{x_{n}}[a_{n}]. \nonumber
\end{eqnarray}
Combining (3.45), (3.46) and direct computations, we obtain
\begin{eqnarray}
&&\sum_{|\alpha|=1}
\text{trace}\Big[\partial_{\xi_{n}}^{2}\pi_{\xi_{n}}^{+}q_{-1}\times
\partial_{x_{n}}\Big(\frac{c(X)}{|\xi|^{2}}-\frac{2g(X,\xi)c(\xi)}{|\xi|^{4}}\Big)
\Big]\nonumber\\
&=&\frac{1} {(\xi_{n}-i)^{3}(1+\xi_{n}^{2})}{\rm tr}\Big[c(\xi')\partial_{x_{n}}[c(X)]\Big]
-\frac{h'(0)} {(\xi_{n}-i)^{3}(1+\xi_{n}^{2})^{2}}{\rm tr}\Big[c(\xi')c(X)\Big]\nonumber\\
&&+\frac{8} {(\xi_{n}-i)^{3}(1+\xi_{n}^{2})^{2}} \partial_{x_{n}}\Big(g(X,\xi')\Big)
+\frac{4\xi_{n}^{3}-8} {(\xi_{n}-i)^{3}(1+\xi_{n}^{2})^{3}} h'(0)g(X,\xi')\nonumber\\
&&+\frac{8\xi_{n}} {(\xi_{n}-i)^{5}(\xi_{n}+i)^{2}} \partial_{x_{n}}(a_{n})
+\frac{4\xi_{n}^{3}-8\xi_{n}} {(\xi_{n}-i)^{6}(\xi_{n}+i)^{3}} a_{n}h'(0).
\end{eqnarray}
Considering for $i<n$,  $\int_{|\xi'|=1}\{\xi_{i_1}\xi_{i_2}\cdots\xi_{i_{2d+1}}\}\sigma(\xi')=0$.
Therefore
\begin{eqnarray}
\text{ Case \ (9) }
&=&\Big(-\frac{367}{128}\big(h'(0)\big)^{2}+\frac{103}{64}h''(0)\Big)\pi\Omega_{3}\texttt{d}x' \nonumber\\
&&+\frac{1}{2}\int_{|\xi'|=1}\int_{-\infty}^{+\infty}\sum_{|\alpha|=1}
\text{trace}\Big[\partial_{\xi_{n}}^{2}\pi_{\xi_{n}}^{+}q_{-1}\times
\partial_{x_{n}}\Big(\frac{c(X)}{|\xi|^{2}}-\frac{2g(X,\xi)c(\xi)}{|\xi|^{4}}\Big)
\Big](x_{0})\texttt{d}\xi_{n}\sigma(\xi')\texttt{d}x'\nonumber\\
&&=\Big(-\frac{367}{128}\big(h'(0)\big)^{2}+\frac{103}{64}h''(0)\Big)\pi\Omega_{3}\texttt{d}x'
-\frac{3}{8}\pi\Omega_{3} \partial_{x_{n}}(a_{n})\texttt{d}x'+\frac{15}{64}\pi a_{n}h'(0)\Omega_{3}\texttt{d}x'.
\end{eqnarray}

\textbf{Case (10)}: \ $r=-2, \ \ell=-1, \ k=0, \ j=1, \ |\alpha|=0$

From (2.9), we have
\begin{equation}
\text{ Case \ (10)}=-\frac{1}{2}\int_{|\xi'|=1}\int_{-\infty}^{+\infty}
\text{trace}\Big[\partial_{x_{n}}\pi_{\xi_{n}}^{+}q_{-2}\times
\partial_{\xi_{n}}^{2}q_{-1}\Big](x_{0})\texttt{d}\xi_{n}\sigma(\xi')\texttt{d}x'.
\end{equation}
By the Leibniz rule, trace property and "++" and "-~-" vanishing
after the integration over $\xi_n$ in \cite{FGLS}, then
\begin{eqnarray}
&&\int^{+\infty}_{-\infty}{\rm trace}
\Big[\partial_{x_{n}}\pi_{\xi_{n}}^{+} q_{-2} \times
\partial_{\xi_{n}}^{2}q_{-1}(D^{-1})\Big]\texttt{d}\xi_n \nonumber\\
&=& \int^{+\infty}_{-\infty}{\rm trace}
\Big[\partial_{x_{n}}q_{-2}(D^{-1})\times\partial_{\xi_{n}}^{2}q_{-1}(D^{-1})\Big]\texttt{d}\xi_n
-\int^{+\infty}_{-\infty}{\rm trace}
\Big[\partial_{x_{n}}q_{-2}(D^{-1})\times\partial_{\xi_{n}}^{2}\pi_{\xi_{n}}^{+}q_{-1}(D^{-1})\Big]\texttt{d}\xi_n.\nonumber\\
\end{eqnarray}
Combining these assertions, we see
\begin{eqnarray}
\text{ Case \ (10)}&=&\text{ Case \ (9)}-\frac{1}{2}\int_{|\xi'|=1}\int^{+\infty}_{-\infty}{\rm trace}
\Big[\partial_{x_{n}}q_{-2} \times\partial_{\xi_{n}}^{2}q_{-1} \Big]\texttt{d}\xi_{n}\sigma(\xi')\texttt{d}x'\nonumber\\
&=&\text{ Case \ (9)}-\frac{1}{2}\int_{|\xi'|=1}\int^{+\infty}_{-\infty}{\rm trace}
\Big[\partial_{x_{n}}\sigma_{-2}(D^{-1})\times \partial_{\xi_{n}}^{2}q_{-1} \Big]\texttt{d}\xi_{n}\sigma(\xi')\texttt{d}x'\nonumber\\
&&-\frac{1}{2}\int_{|\xi'|=1}\int^{+\infty}_{-\infty}{\rm trace}
\Big[\partial_{x_{n}}\Big(\frac{c(X)}{|\xi|^{2}}-\frac{2g(X,\xi)c(\xi)}{|\xi|^{4}}\Big)
\times\partial_{\xi_{n}}^{2}q_{-1} \Big]\texttt{d}\xi_{n}\sigma(\xi')\texttt{d}x'.
\end{eqnarray}
By Lemma 3.2, a simple computation shows
\begin{equation}
\partial_{\xi_{n}}^{2}q_{-1}(x_{0})\Big|_{|\xi'|=1}
=\frac{6\xi_{n}^{2}-2}{(1+\xi_{n}^{2})^{3}}\sqrt{-1}c(\xi')+\frac{2\xi_{n}^{3}
-6\xi_{n}}{(1+\xi_{n}^{2})^{3}}\sqrt{-1}c(\texttt{d}x_{n}).
\end{equation}
From Lemma 3.2, a simple computation shows
\begin{eqnarray}
\partial_{x_{n}}\Big(\frac{c(X)}{|\xi|^{2}}-\frac{2g(X,\xi)c(\xi)}{|\xi|^{4}}\Big)
&=&\partial_{x_{n}}\Big(\frac{c(X)}{|\xi|^{2}}\Big)
-\partial_{x_{n}}\Big(\frac{2g(X,\xi)c(\xi)}{|\xi|^{4}}\Big)\nonumber\\
&=&\frac{1}{|\xi|^{2}}\partial_{x_{n}}\Big(c(X)\Big)-\frac{h'(0)|\xi'|^{2}}{|\xi|^{4}}c(X)
-\frac{2c(\xi)}{|\xi|^{4}}\partial_{x_{n}}\Big(g(X,\xi)\Big)\nonumber\\
&&-\frac{2g(X,\xi)}{|\xi|^{4}}\partial_{x_{n}}\Big(c(\xi)\Big)
+\frac{4g(X,\xi)c(\xi)}{|\xi|^{6}}\partial_{x_{n}}\Big(|\xi|^{2})\Big).
\end{eqnarray}
Combining (3.52), (3.53) and direct computations, we obtain
\begin{eqnarray}
&&{\rm trace}
\Big[\partial_{x_{n}}\Big(\frac{c(X)}{|\xi|^{2}}-\frac{2g(X,\xi)c(\xi)}{|\xi|^{4}}\Big)
\times\partial_{\xi_{n}}^{2}q_{-1} \Big]\nonumber\\
&=&\frac{6 i\xi_{n}^{2}-2i} {(1+\xi_{n}^{2})^{4}}{\rm tr}\Big[c(\xi')\partial_{x_{n}}[c(X)]\Big]
-\frac{(6 i\xi_{n}^{2}-2i) h'(0)} {(1+\xi_{n}^{2})^{5}}{\rm tr}\Big[c(\xi')c(X)\Big]\nonumber\\
&&+\frac{8(6 i\xi_{n}^{2}-2i)} {(1+\xi_{n}^{2})^{5}} \partial_{x_{n}}\Big(g(X,\xi')\Big)
+\frac{8\xi_{n}(2i\xi_{n}^{3}-6i\xi_{n})} {(1+\xi_{n}^{2})^{5}} \partial_{x_{n}}\Big(g(X,\xi')\Big)\nonumber\\
&&+\frac{-8 i \xi_{n}^{4}-224i\xi_{n}^{2}+24i} {(1+\xi_{n}^{2})^{6}} h'(0)g(X,\xi')
+\frac{8i(\xi_{n}^{5}+2\xi_{n}^{3}+\xi_{n})} {(\xi_{n}-i)^{5}(\xi_{n}+i)^{5}} \partial_{x_{n}}(a_{n})\nonumber\\
&&+ \Big(\frac{8i(-\xi_{n}^{5}-28\xi_{n}^{3}+3 \xi_{n})} {(\xi_{n}-i)^{6}(\xi_{n}+i)^{6}}
+\frac{8i(\xi_{n}^{3}-3\xi_{n})} {(\xi_{n}-i)^{5}(\xi_{n}+i)^{5}}  \Big) a_{n}h'(0).
\end{eqnarray}
Considering for $i<n$,  $\int_{|\xi'|=1}\{\xi_{i_1}\xi_{i_2}\cdots\xi_{i_{2d+1}}\}\sigma(\xi')=0$.
Similar to Case (9), we obtain
\begin{equation}
\text{ Case \ (10) }\Big(-\frac{367}{128}\big(h'(0)\big)^{2}+\frac{103}{64}h''(0)\Big)\pi\Omega_{3}\texttt{d}x'
-\frac{3}{8}\pi\Omega_{3} \partial_{x_{n}}(a_{n})\texttt{d}x'+\frac{15}{64}\pi a_{n}h'(0)\Omega_{3}\texttt{d}x'.
\end{equation}

\textbf{Case (11)}: \ $r=-2, \ \ell=-1, \ k=0, \ j=0, \ |\alpha|=1$

From (2.9), we have
\begin{equation}
\text{ Case \ (11)}=-\int_{|\xi'|=1}\int_{-\infty}^{+\infty}\sum_{|\alpha|=1}
\text{trace}\Big[\partial_{\xi'}^{\alpha}\pi_{\xi_{n}}^{+}q_{-2}\times
\partial_{x'}^{\alpha}\partial_{\xi_{n}}q_{-1}\Big](x_{0})\texttt{d}\xi_{n}\sigma(\xi')\texttt{d}x' .
\end{equation}
By Lemma3.2 and Lemma 3.3, for $i<n$, we have
\begin{equation}
\partial_{x_i}q_{-1}(x_0)=\partial_{x_i}\left(\frac{\sqrt{-1}c(\xi)}{|\xi|^2}\right)(x_0)=
\frac{\sqrt{-1}\partial_{x_i}[c(\xi)](x_0)}{|\xi|^2}
-\frac{\sqrt{-1}c(\xi)\partial_{x_i}(|\xi|^2)(x_0)}{|\xi|^4}=0.
\end{equation}
So Case (11) vanishes.

\textbf{Case (12)}: \ $r=-2, \ \ell=-1, \ k=1, \ j=0, \ |\alpha|=0$

From (2.9) and the Leibniz rule, trace property and "++" and "-~-" vanishing
after the integration over $\xi_n$ in \cite{FGLS}, we have
\begin{eqnarray}
\text{ Case \ (12)}&=&-\frac{1}{2}\int_{|\xi'|=1}\int_{-\infty}^{+\infty}
\text{trace}\Big[\partial_{\xi_{n}}\pi_{\xi_{n}}^{+}q_{-2}\times
\partial_{\xi_{n}}\partial_{x_{n}}q_{-1}\Big](x_{0})\texttt{d}\xi_{n}\sigma(\xi')\texttt{d}x'\nonumber\\
&=&\frac{1}{2}\int_{|\xi'|=1}\int_{-\infty}^{+\infty}
\text{trace}\Big[\pi_{\xi_{n}}^{+}q_{-2}\times
\partial_{\xi_{n}}^{2}\partial_{x_{n}}q_{-1}\Big](x_{0})\texttt{d}\xi_{n}\sigma(\xi')\texttt{d}x'\nonumber\\
&=&\text{ Case \ (7)}+\frac{1}{2}\int_{|\xi'|=1}\int_{-\infty}^{+\infty}
\text{trace}\Big[q_{-2}\times\partial_{\xi_{n}}^{2}\partial_{x_{n}}q_{-1} \Big](x_{0})\texttt{d}\xi_{n}
\sigma(\xi')\texttt{d}x'\nonumber\\
&=&\text{ Case \ (7)}+\frac{1}{2}\int_{|\xi'|=1}\int_{-\infty}^{+\infty}
\text{trace}\Big[\sigma_{-2}(D^{-1})\times\partial_{\xi_{n}}^{2}\partial_{x_{n}}q_{-1} \Big](x_{0})\texttt{d}\xi_{n}
\sigma(\xi')\texttt{d}x'\nonumber\\
&&+\frac{1}{2}\int_{|\xi'|=1}\int_{-\infty}^{+\infty}
\text{trace}\Big[\Big(\frac{c(X)}{|\xi|^{2}}-\frac{2g(X,\xi)c(\xi)}{|\xi|^{4}}\Big) \times
  \partial_{\xi_{n}}^{2}\partial_{x_{n}}q_{-1} \Big](x_{0})\texttt{d}\xi_{n}\sigma(\xi')\texttt{d}x'.\nonumber\\
\end{eqnarray}
From Lemma 3.2-Lemma 3.4 and direct computations, we obtain
\begin{eqnarray}
\partial_{\xi_{n}}^{2}\partial_{x_{n}}q_{-1}(x_0)|_{|\xi'|=1}
&=&\frac{6i\xi_n^{2}-2i}{(1+\xi_n^{2})^{3}}\partial_{x_n}[c(\xi')](x_0)
+\sqrt{-1}h'(0)\Big[\frac{4(1-5\xi_n^{2})}{(1+\xi_n^{2})^{4}}c(\xi')\nonumber\\
&&-\frac{12\xi_n(\xi_n^{2}-1)}{(1+\xi_n^{2})^{4}}c(\texttt{d}x_n)\Big].
\end{eqnarray}
Combining (3.58), (3.59) and direct computations, we obtain
\begin{eqnarray}
&&\text{trace}\Big[\Big(\frac{c(X)}{|\xi|^{2}}-\frac{2g(X,\xi)c(\xi)}{|\xi|^{4}}\Big) \times
  \partial_{\xi_{n}}^{2}\partial_{x_{n}}q_{-1} \Big]\nonumber\\
&=&\frac{(6i\xi_{n}^{2}-2i)} {(1+\xi_{n}^{2})^{4}}{\rm tr}\Big[c(X)\partial_{x_{n}}[c(\xi')]\Big]
+\frac{4i(1-5\xi_{n}^{2})} {(1+\xi_{n}^{2})^{5}}h'(0){\rm tr}\Big[c(\xi')c(X)\Big]
+ \frac{8i(\xi_{n}^{5}-10\xi_{n}^{3}-3 \xi_{n})} {(1+\xi_{n}^{2})^{6}}
 a_{n}h'(0).\nonumber\\
\end{eqnarray}
Therefore
\begin{equation}
\text{ Case \ (12)}=\frac{39}{32}\big(h'(0)\big)^{2}\pi\Omega_{3}\texttt{d}x'-\frac{5}{8}a_{n}h'(0)\pi\Omega_{3}\texttt{d}x'.
\end{equation}

\textbf{Case (13)}: \ $r=-2, \ \ell=-2, \ k=0, \ j=0, \ |\alpha|=0$

From (2.9) and the Leibniz rule , we have
\begin{eqnarray}
\text{ Case \ (13)}&=&-i\int_{|\xi'|=1}\int_{-\infty}^{+\infty}
\text{trace}\Big[\pi_{\xi_{n}}^{+}q_{-2}\times
\partial_{\xi_{n}}q_{-2}\Big](x_{0})\texttt{d}\xi_{n}\sigma(\xi')\texttt{d}x'\nonumber\\
&=&i\int_{|\xi'|=1}\int_{-\infty}^{+\infty}
\text{trace}\Big[\partial_{\xi_{n}}\pi_{\xi_{n}}^{+}\sigma_{-2}(D^{-1})\times
\sigma_{-2}(D^{-1})\Big](x_{0})\texttt{d}\xi_{n}\sigma(\xi')\texttt{d}x'\nonumber\\
&&+i\int_{|\xi'|=1}\int_{-\infty}^{+\infty}
\text{trace}\Big[\partial_{\xi_{n}}\pi_{\xi_{n}}^{+} (\sigma_{-2}(D^{-1}))
\times\Big(\frac{c(X)}{|\xi|^{2}}-\frac{2g(X,\xi)c(\xi)}{|\xi|^{4}}\Big)\Big](x_{0})\texttt{d}\xi_{n}\sigma(\xi')\texttt{d}x'\nonumber\\
&&+i\int_{|\xi'|=1}\int_{-\infty}^{+\infty}
\text{trace}\Big[\partial_{\xi_{n}}\pi_{\xi_{n}}^{+} \Big(\frac{c(X)}{|\xi|^{2}}-\frac{2g(X,\xi)c(\xi)}{|\xi|^{4}}\Big)
\times q_{-2}\Big](x_{0})\texttt{d}\xi_{n}\sigma(\xi')\texttt{d}x'.
\end{eqnarray}
By (2.1) in \cite{Wa1} and the Cauchy integral formula, then
\begin{eqnarray}
\pi^+_{\xi_n}\left(\frac{c(\xi)}{|\xi|^4}\right)(x_0)\Big|_{|\xi'|=1}&=&\pi^+_{\xi_n}\left[\frac{c(\xi')+\xi_nc(dx_n)}{(1+\xi_n^2)^2}\right]
=\frac{1}{2\pi i} \lim_{u\rightarrow 0^-}\int_{\Gamma^+}\frac{\frac{c(\xi')+\eta_nc(dx_n)}{(\eta_n+i)^2(\xi_n+iu-\eta_n)}}
{(\eta_n-i)^2}d\eta_n  \nonumber\\
&=&\left[\frac{c(\xi')+\eta_nc(dx_n)}{(\eta_n+i)^2(\xi_n-\eta_n)}\right]^{(1)}|_{\eta_n=i}
=-\frac{ic(\xi')}{4(\xi_n-i)}-\frac{c(\xi')+ic(dx_n)}{4(\xi_n-i)^2}.\\
\pi^+_{\xi_n}\Big(\frac{g(X,\xi)c(\xi)}{|\xi|^{4}}\Big)(x_0)\Big|_{|\xi'|=1}
&=&g(X',\xi')\pi^+_{\xi_n}\Big(\frac{c(\xi)}{|\xi|^{4}}\Big)(x_0)\Big|_{|\xi'|=1}
+a_{n}\pi^+_{\xi_n}\Big(\frac{\xi_{n}c(\xi)}{|\xi|^{4}}\Big)(x_0)\Big|_{|\xi'|=1}\nonumber\\
&=&g(X',\xi')\frac{-i \xi_n-2 }{4(\xi_{n}-i)^{2}}c(\xi')+g(X',\xi')\frac{-i }{4(\xi_{n}-i)^{2}}c(dx_{n})\nonumber\\
&&+\frac{-a_{n}i }{4(\xi_{n}-i)^{2}}c(\xi')+\frac{-a_{n}i\xi_{n} }{4(\xi_{n}-i)^{2}}c(dx_{n}).\\
\pi^+_{\xi_n}\Big(\frac{c(X)}{|\xi|^{2}}\Big)(x_0)\Big|_{|\xi'|=1}
&=&c(X)\pi^+_{\xi_n}\Big(\frac{1}{|\xi|^{2}}\Big)(x_0)\Big|_{|\xi'|=1}=\frac{-i }{2(\xi_{n}-i)}c(X).
\end{eqnarray}
Hence in this case,
\begin{eqnarray}
\partial_{\xi_{n}}\pi^+_{\xi_n}(\sigma_{-2}(D^{-1}))(x_0)\Big|_{|\xi'|=1}&=&h'(0)\frac{-i\xi_n-3}{4(\xi_n-i)^{3}}c(\xi')c(\texttt{d}x_n)c(\xi')
+h'(0)\frac{i}{(\xi_n-i)^{3}}c(\xi')  \nonumber\\
&&+h'(0)\frac{i\xi_n-1}{4(\xi_n-i)^{3}}c(\texttt{d}x_n)+\frac{-i}{2(\xi_n-i)^{3}}\partial_{x_n}[c(\xi')](x_0) \nonumber\\
&&+\frac{i\xi_n+3}{4(\xi_n-i)^{3}} c(\xi') c(\texttt{d}x_n)\partial_{x_n}[c(\xi')](x_0)\nonumber\\
&&+h'(0)\frac{-2i\xi_n-8}{8(\xi_n-i)^{4}}c(\xi')  +h'(0)\frac{i\xi_n^{2}+4\xi_n-9i}{8(\xi_n-i)^{4}}c(\texttt{d}x_n).
\end{eqnarray}
Then
\begin{equation}
i\int_{|\xi'|=1}\int_{-\infty}^{+\infty}
\text{trace}\Big[\partial_{\xi_{n}}\pi_{\xi_{n}}^{+} (\sigma_{-2}(D^{-1}))
\times\Big(\frac{c(X)}{|\xi|^{2}}-\frac{2g(X,\xi)c(\xi)}{|\xi|^{4}}\Big)\Big](x_{0})\texttt{d}\xi_{n}\sigma(\xi')\texttt{d}x'
=\frac{15}{16}\pi a_{n}h'(0)\Omega_{3}\texttt{d}x'.
\end{equation}
Similarly, we obtain
\begin{eqnarray}
&&i\int_{|\xi'|=1}\int_{-\infty}^{+\infty}\text{trace}\Big[\partial_{\xi_{n}}\pi_{\xi_{n}}^{+}
\Big(\frac{c(X)}{|\xi|^{2}}-\frac{2g(X,\xi)c(\xi)}{|\xi|^{4}}\Big)\times q_{-2}\Big](x_{0})\texttt{d}\xi_{n}\sigma(\xi')\texttt{d}x'\nonumber\\
&=&\frac{\pi}{2}|X|^{2}_{g^{TM}}\Omega_{3}\texttt{d}x'+\frac{35\pi}{64}|X'|^{2}_{g^{\partial M}}h'(0)\Omega_{3}\texttt{d}x'
-\pi a_{n}^{2}\Omega_{3}\texttt{d}x'.
\end{eqnarray}
Therefore
\begin{equation}
\text{ Case \ (13)}=-\frac{821}{256}\big(h'(0)\big)^{2}\pi\Omega_{3}\texttt{d}x'+\frac{15}{16}\pi a_{n}h'(0)\Omega_{3}\texttt{d}x'
+\frac{\pi}{2}|X|^{2}_{g^{TM}}\Omega_{3}\texttt{d}x'+\frac{35\pi}{64}|X'|^{2}_{g^{\partial M}}h'(0)\Omega_{3}\texttt{d}x'
-\pi a_{n}^{2}\Omega_{3}\texttt{d}x'.
\end{equation}

\textbf{Case (14)}: \ $r=-1, \ \ell=-3, \ k=0, \ j=0, \ |\alpha|=0$

From (2.9) and the Leibniz rule , we have
\begin{eqnarray}
\text{ Case \ (14)}&=&-i\int_{|\xi'|=1}\int_{-\infty}^{+\infty}
\text{trace}\Big[\pi_{\xi_{n}}^{+}q_{-1}\times
\partial_{\xi_{n}}q_{-3}\Big](x_{0})\texttt{d}\xi_{n}\sigma(\xi')\texttt{d}x'\nonumber\\
&=&i\int_{|\xi'|=1}\int_{-\infty}^{+\infty}
\text{trace}\Big[\partial_{\xi_{n}}\pi_{\xi_{n}}^{+}q_{-1}\times
\sigma_{-3}(D^{-1})\Big](x_{0})\texttt{d}\xi_{n}\sigma(\xi')\texttt{d}x'\nonumber\\
&&+i\int_{|\xi'|=1}\int_{-\infty}^{+\infty}
\text{trace}\Big[\partial_{\xi_{n}}\pi_{\xi_{n}}^{+}q_{-1}\times
R_{-3}\Big|_{|\xi'|=1}\Big](x_{0})\texttt{d}\xi_{n}\sigma(\xi')\texttt{d}x'.
\end{eqnarray}
By (2.2.29) in \cite{Wa3}, we have
\begin{equation}
\partial_{\xi_n}\pi^+_{\xi_n}q_{-1}(x_0)|_{|\xi'|=1}=-\frac{c(\xi')+ic(\texttt{d}x_n)}{2(\xi_n-i)^2}.
\end{equation}
In the orthonormal frame field, we have
\begin{eqnarray}
R_{-3}(x_{0})\Big|_{|\xi'|=1}&=&\frac{-i(\xi_{n}^{4}+\xi_{n}^{2}-2) }{(1+\xi_{n}^{2})^{4}}h'(0)c(\xi)c(X)c(dx_{n})
                                -\frac{i(2\xi_{n}^{3}+4\xi_{n})  }{(1+\xi_{n}^{2})^{4}}h'(0)c(\xi)c(X)c(\xi')\nonumber\\
                             && -\frac{i}{(1+\xi_{n}^{2})^{3}}c(\xi)c(X)c(\xi')c(dx_{n})\partial_{x_{n}}c(\xi')
                             +\frac{i\xi_{n}}{(1+\xi_{n}^{2})^{3}}c(\xi)c(X)\partial_{x_{n}}c(\xi')\nonumber\\
                             &&-\frac{i}{(1+\xi_{n}^{2})^{2}}c(\xi)c(X)
                             +\frac{2i}{(1+\xi_{n}^{2})^{3}}g(X',\xi')c(\xi)c(\xi)
                             +\frac{2ia_{n}\xi_{n}}{(1+\xi_{n}^{2})^{3}}c(\xi)c(\xi)\nonumber\\
                             &&-\frac{i}{(1+\xi_{n}^{2})^{2}} c(\xi)\sum _{j=1}^{n} c(dx_{j})\partial_{x_{j}}c(X)
                             + \frac{ 2i}{(1+\xi_{n}^{2})^{3}}\partial_{x_{j}}(g(X',\xi')) c(\xi)\sum _{j=1}^{n} c(dx_{j})c(\xi) \nonumber\\
                             &&+\frac{2i\xi_{n}\partial_{x_{j}}(a_{n})}{(1+\xi_{n}^{2})^{3}} c(\xi)\sum _{j=1}^{n}c(dx_{j})c(\xi)
                             -\frac{i}{(1+\xi_{n}^{2})^{2}} c(\xi) c(dx_{n})\partial_{x_{n}}c(X)
                            \nonumber\\
                             && +\frac{i h'(0)}{(1+\xi_{n}^{2})^{3}} c(\xi) c(dx_{n})c(X)
                             +\frac{2i}{(1+\xi_{n}^{2})^{3}}\partial_{x_{n}}g(X,\xi) c(\xi) c(dx_{n})c(\xi)
                             \nonumber\\
                             &&+\frac{2i}{(1+\xi_{n}^{2})^{3}} g(X,\xi) c(\xi) c(dx_{n}) \partial_{x_{n}}c(\xi')
                             -\frac{4i h'(0)}{(1+\xi_{n}^{2})^{4}} g(X,\xi) c(\xi) c(dx_{n}) c(\xi).
\end{eqnarray}
By the relation of the Clifford action and ${\rm tr}{AB}={\rm tr }{BA}$, then we have the equalities:
\begin{eqnarray}
&&{\rm tr}\Big[c(\xi')c(\xi)c(X)c(\xi')c(dx_{n})\partial_{x_{n}}c(\xi')\Big]\nonumber\\
&=& {\rm tr}\Big[c(\xi')\Big(c(\xi')+\xi_{n}c(dx_{n}))\Big)c(X)c(\xi')c(dx_{n})\partial_{x_{n}}c(\xi')\Big]\nonumber\\
&=& {\rm tr}\Big[c(\xi')c(\xi')c(X)c(\xi')c(dx_{n})\partial_{x_{n}}c(\xi')\Big]
+ {\rm tr}\Big[c(\xi')\xi_{n}c(dx_{n})c(X)c(\xi')c(dx_{n})\partial_{x_{n}}c(\xi')\Big]\nonumber\\
&=&-2a_{n}h'(0)-{\rm tr}\Big[\partial_{x_{n}}(c(\xi'))c(X)\Big].
\end{eqnarray}
Considering for $i<n$,  $\int_{|\xi'|=1}\{\xi_{i_1}\xi_{i_2}\cdots\xi_{i_{2d+1}}\}\sigma(\xi')=0$.
 From Lemma 3.5, combining (3.71)-(3.73) and direct computations, we obtain
\begin{eqnarray}
&&i\int_{|\xi'|=1}\int_{-\infty}^{+\infty}
\text{trace}\Big[\partial_{\xi_{n}}\pi_{\xi_{n}}^{+}q_{-1}\times
R_{-3}\Big|_{|\xi'|=1}\Big](x_{0})\texttt{d}\xi_{n}\sigma(\xi')\texttt{d}x'\nonumber\\
&=&-\frac{5}{8}\pi a_{n}h'(0)\Omega_{3}\texttt{d}x' +\frac{3 }{4}\pi \partial_{ x_{n}}(a_{n})\Omega_{3}\texttt{d}x'
+\frac{3}{4}\pi\Omega_{3} C_{1}^{1}(DX^{*})\texttt{d}x'.
\end{eqnarray}
Therefore
\begin{eqnarray}
\text{ Case \ (14)}
&=&\Big(\frac{239}{64}\big(h'(0)\big)^{2}-\frac{27}{16}h''(0)-\frac{11}{192}s_{\partial_{M}}\Big)\pi\Omega_{3}\texttt{d}x'\nonumber\\
&&-\frac{5}{8}\pi a_{n}h'(0)\Omega_{3}\texttt{d}x' +\frac{3 }{4}\pi \partial_{ x_{n}}(a_{n})\Omega_{3}\texttt{d}x'
+\frac{3}{4}\pi\Omega_{3} C_{1}^{1}(DX^{*})\texttt{d}x'.
\end{eqnarray}

\textbf{Case (15)}: \ $r=-3, \ \ell=-1, \ k=0, \ j=0, \ |\alpha|=0$

From (2.9) we have
\begin{equation}
\text{ Case \ (15)}=-i\int_{|\xi'|=1}\int_{-\infty}^{+\infty}
\text{trace}\Big[\pi_{\xi_{n}}^{+}q_{-3}\times
\partial_{\xi_{n}}q_{-1}\Big](x_{0})\texttt{d}\xi_{n}\sigma(\xi')\texttt{d}x'.
\end{equation}
By the Leibniz rule, trace property and "++" and "-~-" vanishing
after the integration over $\xi_n$ in \cite{FGLS}, then
\begin{eqnarray}
&&\int^{+\infty}_{-\infty}{\rm trace}
\Big[\pi_{\xi_{n}}^{+}q_{-3}\times
\partial_{\xi_{n}}q_{-1}\Big]\texttt{d}\xi_n \nonumber\\
&=& \int^{+\infty}_{-\infty}{\rm trace}
\Big[q_{-3}\times
\partial_{\xi_{n}}q_{-1}\Big]\texttt{d}\xi_n
-\int^{+\infty}_{-\infty}{\rm trace}
\Big[q_{-3}\times
\partial_{\xi_{n}}\pi_{\xi_{n}}^{+}q_{-1}\Big]\texttt{d}\xi_n.
\end{eqnarray}
Combining these assertions, we see
\begin{equation}
\text{ Case \ (15)}=\text{ Case \ (14)}-i\int_{|\xi'|=1}\int_{-\infty}^{+\infty}
\text{trace}\Big[R_{-3}\times
\partial_{\xi_{n}}q_{-1}\Big](x_{0})\texttt{d}\xi_{n}\sigma(\xi')\texttt{d}x'.
\end{equation}
By Lemma 3.2, a simple computation shows
\begin{equation}
\partial_{\xi_{n}}q_{-1}(x_{0})\Big|_{|\xi'|=1}
=\frac{-2\xi_{n}}{(1+\xi_{n}^{2})^{2}}\sqrt{-1}c(\xi')+\frac{1-\xi_{n}^{2}}{(1+\xi_{n}^{2})^{2}}\sqrt{-1}c(\texttt{d}x_{n}).
\end{equation}
Similar to Case (14), combining (3.72),(3.79) and direct computations, we obtain
\begin{eqnarray}
\text{ Case \ (15)}
&=&\Big(\frac{239}{64}\big(h'(0)\big)^{2}-\frac{27}{16}h''(0)-\frac{11}{192}s_{\partial_{M}}\Big)\pi\Omega_{3}\texttt{d}x'\nonumber\\
&&-\frac{5}{8}\pi a_{n}h'(0)\Omega_{3}\texttt{d}x' +\frac{3 }{4}\pi \partial_{ x_{n}}(a_{n})\Omega_{3}\texttt{d}x'
+\frac{3}{4}\pi\Omega_{3}  C_{1}^{1}(\nabla^{\partial M} (X'|_{\partial M})^{*})\texttt{d}x'.
\end{eqnarray}
Now $\Phi$  is the sum of the \textbf{case ($1,2,\cdots,15$)} , then we obtain
\begin{eqnarray}
\sum_{I=1}^{15} \textbf{case I}
&=&\Big(\frac{399}{256}\big(h'(0)\big)^{2}-\frac{29}{32}h''(0)-\big(\frac{17}{96}
+\frac{5}{32}i\big)s_{\partial_{M}}\Big)\pi\Omega_{3}\texttt{d}x'\nonumber\\
&&+\frac{25}{32}\pi a_{n}h'(0)\Omega_{3}\texttt{d}x'
-\frac{\pi}{2}|X|^{2}_{g^{TM}}\Omega_{3}\texttt{d}x'+\frac{35\pi}{64}|X'|^{2}_{g^{\partial M}}h'(0)\Omega_{3}\texttt{d}x'
-\pi a_{n}^{2}\Omega_{3}\texttt{d}x'\nonumber\\
&&-\frac{3 }{4}\pi \partial_{ x_{n}}(a_{n})\Omega_{3}\texttt{d}x'
+\frac{15}{16}\pi C_{1}^{1}(\nabla^{\partial M} (X'|_{\partial M})^{*})\Omega_{3}\texttt{d}x'.
\end{eqnarray}
Hence we conclude that,
\begin{thm}
 Let $M$ be a five dimensional compact manifold  with the boundary $\partial M$, and
 the Dirac operators with one-form perturbations $\widetilde{D}=D+c(X)$, then
\begin{eqnarray}
Vol_{5}^{(1, 1)}
&=&\int_{\partial_{M}}\Big[\frac{399}{256}\big(h'(0)\big)^{2}-\frac{29}{32}h''(0)-\big(\frac{17}{96}
+\frac{5}{32}i\big)s_{\partial_{M}}\nonumber\\
&&+\frac{25}{32} a_{n}h'(0)
-\frac{1}{2}|X|^{2}_{g^{TM}}+\frac{35}{64}|X'|^{2}_{g^{\partial M}}h'(0)
- a_{n}^{2}\nonumber\\
&&-\frac{3 }{4} \partial_{ x_{n}}(a_{n})
+\frac{15}{16} C_{1}^{1}(\nabla^{\partial M} (X'|_{\partial M})^{*})\Big]\pi \Omega_{3}{\rm dvol}_{\partial_{M}}.
\end{eqnarray}
\end{thm}

Next we recall the Einstein-Hilbert action for manifolds with boundary  (see \cite{Wa3} or \cite{Wa4}),
\begin{equation}
I_{\rm Gr}=\frac{1}{16\pi}\int_Ms{\rm dvol}_M+2\int_{\partial M}K{\rm dvol}_{\partial_M}:=I_{\rm {Gr,i}}+I_{\rm {Gr,b}},
\end{equation}
  where
  \begin{equation}
K=\sum_{1\leq i,j\leq {n-1}}K_{i,j}g_{\partial M}^{i,j};~~K_{i,j}=-\Gamma^n_{i,j},
\end{equation}
and $K_{i,j}$ is the second fundamental form, or extrinsic
curvature. Take the metric in Section 2, then by Lemma A.2 in \cite{Wa3},
$K_{i,j}(x_0)=-\Gamma^n_{i,j}(x_0)=-\frac{1}{2}h'(0),$ when $i=j<n$,
otherwise is zero. For $n=5$, then
  \begin{equation}
K(x_0)=\sum_{i,j}K_{i.j}(x_0)g_{\partial M}^{i,j}(x_0)=\sum_{i=1}^{4}K_{i,i}(x_0)=-2h'(0).
\end{equation}
Then
   \begin{equation}
I_{\rm {Gr,b}}=-4h'(0){\rm Vol}_{\partial M}.
\end{equation}
On the other hand, by Proposition 2.10 in \cite{Wa5} and Th 3.9 in \cite{WW1}, we have
\begin{lem}\cite{WW1}
 Let M be a $5$-dimensional compact manifold  with the boundary $\partial M$, then
 \begin{equation}
s_{M}(x_{0})=3\big(h'(0)\big)^{2}-4h''(0)+s_{\partial_{M}}(x_{0}).
\end{equation}
\end{lem}
Hence from (3.82)-(3.87), we obtain
\begin{thm}
 Let $M$ be a five dimensional compact spin manifold  with the boundary $\partial M$, and
 the Dirac operators with one-form perturbations $\widetilde{D}=D+c(X)$.
  Let $X=X'+a_{n}dx_{n}$ near the boundary and  $X'$ is a one form on $\partial_{M}$,
 the following identity holds
\begin{eqnarray}
 \widetilde{{\rm Wres}}[\pi^+(D+c(X))^{-1}\circ \pi^+(D+c(X))^{-1}]
&=&\int_{\partial_{M}}\Big[
\frac{1}{16} \Big( \frac{225}{32}K^{2}+\frac{29}{4}s_{M}\big|_{\partial_{M}}-\big(\frac{155}{12}+5i \big)s_{\partial_{M}}\Big)\nonumber\\
&&-\frac{25}{32} a_{n}K
- |X'|^{2}_{g^{\partial M}}-\frac{35}{64}|X'|^{2}_{g^{\partial M}}K
- 3a_{n}^{2}|_{\partial M}\nonumber\\
&&-\frac{3 }{2} \partial_{ x_{n}}(a_{n})|_{\partial M}
+\frac{15}{8} C_{1}^{1}(\nabla^{\partial M} (X'|_{\partial M})^{*})\Big]\pi^{3}{\rm dvol}_{\partial_{M}},
\end{eqnarray}
where $s_{M}$, $s_{\partial_{M}}$ are respectively scalar curvatures on $M$ and $\partial_{M}$,
 and the vector field $(X'|_{\partial M})^{*}$ is the metric dual of $X'|_{\partial M}$,
  $\nabla^{\partial M}$ is the Levi-civita connection on $\partial M$, $C_{1}^{1}$ is
  the contraction of $(1,1)$ tensors.
\end{thm}

\section*{ Acknowledgements}
~The work of the first author was supported by NSFC. 11501414. The work of the second author was supported by NSFC. 11771070.
 The authors also thank the referee
for his (or her) careful reading and helpful comments.

\end{document}